%% file: preprint.tex
\documentclass[
hidempi,
hidelinks
]{mpi2015-cscpreprint}

\input{macros}

\begin{document}

\title{%
  \texorpdfstring{Reduced rank extrapolation\\ for multi-term Sylvester equations}
                 {Reduced rank extrapolation for multi-term Sylvester equations}
}

\author[$\ast$]{Peter Benner~\orcidlink{0000-0003-3362-4103}} % chktex 8
\author[$\dagger$]{Pascal~den~Boef~\orcidlink{0009-0004-9061-9284}} % chktex 8
\author[$\ddagger$]{Patrick~K\"urschner~\orcidlink{0000-0002-6114-8821}} % chktex 8
\author[$\ast$]{Xiaobo~Liu~\orcidlink{0000-0001-8470-8388}} % chktex 8
\author[$\ast$]{Jens~Saak~\orcidlink{0000-0001-5567-9637}} % chktex 8

\affil[$\ast$]{Max Planck Institute for Dynamics of Complex Technical Systems, Magdeburg, Germany.}
\affil[$\dagger$]{Eindhoven University of Technology, Eindhoven, Netherlands.}
\affil[$\ddagger$]{Leipzig University of Applied Sciences, Leipzig, Germany.}

\abstract{\input{abstract}}



\novelty{\input{novelty}}

\maketitle

\input{main.tex}

%\vspace{0.5cm}
%\textbf{Acknowledgment} $\cdot$

%\appendix
%\input{appendix.tex}

\bibliographystyle{plainurl}
\bibliography{journals, reference}

\end{document}

%% file: abstract.tex
We investigate the acceleration of stationary iterations for multi-term Sylvester equation by means of reduced rank extrapolation (RRE).
Theoretical convergence results and implementations are provided for both small and large-scale problems. 
For the large-scale problems, an inexact non-stationary iteration is discussed, which makes use of low-rank matrix approximations.
Numerical experiments illustrate the potential of the RRE acceleration which often leads to a substantial gain in convergence speed and therefore reducing the consumption of storage and computing time.

%% file: novelty.tex
An RRE-accelerated framework for multi-term Sylvester and Lyapunov-plus-positive equations is developed, providing:
(i) theoretical convergence results for iterative fixed-point methods; and
(ii) an efficient implementation for inexact, non-stationary fixed-point processes, combining RRE with dynamically adjusted inexact Sylvester solves and practical rank-truncation techniques.

%% file: main.tex
% chktex-file 3  % disable {} warnings for exponents etc.
% chktex-file 17 % disable warning about unmatched numbers of opening and
% closing brackets
\section{Introduction}

The general linear matrix equation
\begin{align}\label{genmateq}
  \cA(X)=-Y
\end{align}
is defined by a linear, matrix valued operator $\cA:\Rnm\mapsto \Rnm$.
We investigate the case when $\cA$ can be decomposed additively as
\begin{align}\label{sylmateq}
  \cA(X)=\cL(X)+\Pi(X).
\end{align}
The $\cL(X)$ represents a Sylvester operator
\begin{align}\label{standardsylv}
  \cL(X)=AX+XB,
\end{align}
and the second operator is given by
\begin{equation}\label{posoper}
  \Pi(X)=\sum\limits_{k=1}^{\ell}N_{k}XH_{k}, \quad
  N_k\in\Rnn,~H_k\in\Rmm,
\end{equation}
so equation~\eqref{genmateq} represents a generalized (multi-term) Sylvester equation.
Note that it is possible to define the Sylvester operators in a more general form as
\begin{align}\label{gensylvop}
  \cL(X)=AXC+MXB,
\end{align}
where $C\in\Rmm$ and $M\in\Rnn$ are invertible.
For the sake of simpler notations, we mostly restrict to the case \(C=I\) and \(M=I\), but the investigation carried out here also extends to~\eqref{gensylvop}.

If $B=A^{\TT}$, then $\cL(X)$ becomes a Lyapunov operator (additionally $C=M^{\TT}$ for~\eqref{gensylvop}). Moreover, with $H_k=N_k^{\TT}$, $Y=Y^{\TT}$, equation~\eqref{genmateq} is a generalized Lyapunov equation admitting a symmetric solution. In control theory applications, this is frequently referred to as a \emph{Lyapunov-plus-positive} equation~\cite{BenS13}.

There exist established numerical methods for solving standard Lyapunov and Sylvester matrix equations $\cL(X)=-Y$.
In the small and dense case, the Schur decomposition method~\cite{BarS72} and its Hessenberg variant~\cite{GolNL79}
are the standard approaches. For medium- to large-scale problems, iterative
solvers are in general preferred due to the high computational and memory
overhead of factorization-based methods, and often, the low numerical rank of
the solution is exploited to obtain a factored solution; see, e.g.,
\cite{BenS13, Sim16a} and references therein.% chktex 2

The Lyapunov-plus-positive variant of~\eqref{genmateq} plays an important role in control and model reduction of bilinear and special linear stochastic differential equations~\cite{BenD11, Kle69}. Different numerical methods have been designed for solving~\eqref{genmateq},
including factorization-based methods which exploit the Sherman-Morrison-Woodbury formula~\cite{ColDR93,Dam08,HaoS21}, splitting-based iterations~\cite{Dam08, ShaSS15}, and other low-rank schemes~\cite{BenB13}. On the other hand, Krylov-based methods have been developed
in, e.g.,~\cite{KreT11,JaMPR18,PalK2021,PalCS25}.

We are mainly interested in the situation when $\cL$ and $\Pi$ in
\eqref{genmateq} represent a convergent splitting of $\cA$, e.g., of the form % chktex 2
that the spectral radius bound $\rho(\cL^{-1}\Pi)<1$ holds. More stringently, the requirement can become $\|\cL^{-1}\Pi\|<1$ for some matrix norm. This condition is necessary for the solution $X$ to be positive semidefinite in the Lyapunov-plus-positive case~\cite{BenB13}.
As a consequence, we then have that
\[
  X=\cL^{-1}(-Y-\Pi(X))
\]
is a fixed-point equation and admits the application of basic and splitting-based iterations
\begin{align}\label{eq:basicit}
  X_{k+1}=\cL^{-1}(-Y-\Pi(X_{k})),\quad k\geq 0.
\end{align}
This splitting-based fixed-point iteration has been studied for the
Lyapunov-plus-positive case in, e.g.,~\cite{Dam08, ShaSS15}.
For the general case, the invertibility of $\cL$ is ensured by the spectral
condition $\Lambda(A)\cap\Lambda(-B)=\emptyset$; see~\cite[Thm.~4.4.6]{HorJ91}.

In particular, this condition holds if $\Lambda(A), \Lambda(B)\subset \C_-$ (or $\C_+$).
Hence, in each step of the iteration~\eqref{eq:basicit}
one has to solve a standard Sylvester equation
\begin{equation}\label{eq:basicsylv}
  AX_{k+1}+X_{k+1}B=-Y-\Pi(X_{k}).
\end{equation}

The method of reduced rank extrapolation (RRE) was independently proposed
in~\cite{eddy1979extrapolating, kaniel1974least, mevsina1977convergence}, in the
1970s. It was originally introduced as a vector extrapolation method for accelerating the solutions of linear system of equations,
in the context of finding the solution of a fixed-point iteration
\begin{equation}\label{eq:fixpoint}
  f\colon \mathbb{R}^n \to \mathbb{R}^n,\quad
  x_{k+1} = f(x_k), \ k=0,1,\dots,
\end{equation}
for some initial vector $x_0\in\mathbb{R}^n$.
Let
\[\Gamma:=\left\lbrace
    \gamma\in\R^{w}\ \middle|\ \sum_{i=1}^{w} \gamma_i = 1 \right\rbrace.
\]
The RRE constructs an extrapolant vector $\widehat{x}$ as a linear combination of $w$ (which is called the \emph{window size}) iterates, say, $\{ x_i \}_{i=1}^{w}$:
\begin{equation}\label{eq:rre-extrap}
  \widehat{x} = \sum_{i=1}^{w} \widehat{\gamma_i} x_i,
\end{equation}
where the $\widehat{\gamma_i}$ are the solution to the constrained optimization
problem
\begin{equation*}
  \widehat{\gamma} = \underset{\gamma\in\Gamma}{\argmin}
  \left\| \sum_{i=1}^{w} \gamma_i (x_{i+1}-x_i) \right\|_2.
\end{equation*}
Let the matrix of first-order differences of the $x_i$ be defined as
\begin{equation}\label{eq:def-ui}
  U_w = [u_1, u_2, \ldots, u_{w}], \quad u_i:= x_{i+1}-x_i.
\end{equation}
Then, this minimization problem can be expressed compactly as
\begin{equation}\label{eq:rre-minprob}
  \hgamma = \underset{\gamma\in\Gamma}{\argmin}
  \left\| U_w\gamma \right\|_2.
\end{equation}
Note that $w+1$ successive iterates are required for calculating a
$w$-term extrapolant.
The optimization problem~\eqref{eq:rre-minprob} is convex and can be transformed into a (unconstrained) least-squares problem, which can be solved via QR decomposition~\cite{sidi1991efficient}, see also Section~\ref{ssec:conv_small}.
We refer to~\cite{Saad25} for a survey on RRE and similar acceleration techniques, also for nonlinear and other problems.

Several algorithmic modes of RRE exist, notably cycling and non-cycling (see, e.g.,~\cite{sidi2020convergence} in which these modes are referred to as n-Mode and C-mode, respectively).
In the cycling mode, which we consider in this work, a cycle refers to the collection of enough iterates of the fixed-point process, computing the extrapolant $\hx$, and restarting with initial guess $x_0=\hx$.
(The non-cycling mode omits the restarting and hence can be seen as fully non-intrusive mode applicable for problems with limited access to the fixed-point process.)
Finally, we note that RRE has been generalized to accelerate the solution of
matrix-valued fixed-point processes for standard Lyapunov and algebraic Riccati
equations, where the iterates are symmetric low-rank factorized matrices in~\cite{DeKLM25}.
In that work, the arising fixed-point type iterations are \emph{stationary} (as~\eqref{eq:fixpoint}) or \emph{non-stationary} (as~\eqref{eq:fixpoint} with iteration-dependent $f$).

In this paper, we investigate the application of RRE for the numerical solution of generalized Sylvester equations. More specifically, our contributions are as follows:
\begin{enumerate}
\item We propose a novel extension of RRE to accelerate the above iteration~\eqref{eq:basicit} for, both, small- and large-scale generalized Sylvester equations.
\item We derive several theoretical convergence results for the iteration~\eqref{eq:basicit} with or without RRE\@. For applicability to large-scale problems, both exact and inexact solves are considered in the analysis, leading to stationary and non-stationary fixed-point processes, respectively.
\item We formulate algorithms for the proposed methods including practical considerations such as rank truncation and residual norm estimation.
\end{enumerate}

Throughout the paper, a tilde is used to indicate approximate quantities, while extrapolated quantities wear a hat.
Our analysis and implementations focus on equations containing real matrices; however, they can be straightforwardly generalized to the complex case. In particular, the RRE is formulated in a real vector space, as the extrapolation coefficients and the extrapolated vector remain real when the underlying matrix equation involves only real matrices.

%%%%%%%%%%%%%%%%%%%%%%%%%%%%
\section{Analysis and algorithm for small, dense problems}
\begin{algorithm2e}[tbp]
  \caption{%
    Stationary iteration~\eqref{eq:basicit} with RRE for the generalized
    Sylvester equation~\eqref{sylmateq}.  }%
  \label{alg:basic}
  \KwIn{%
    Matrices $A,~B,~Y,~\{N_i\}_{i=1}^{\ell},~\{H_i\}_{i=1}^{\ell},~X_0$.  }
  \KwOut{%
    Approximate solution~$\tX\approx X$.  }

  $k=1$\; \While{not converged}{ Solve
    $AX_{k}+X_{k}B=-Y-\sum\limits_{i=1}^{\ell}N_i X_{k-1} H_i$\;
    \If{$k\geq w~\land~\mod(k,w)=0$}{ Extrapolation step: $X_{k}$ =
      RRE$_{w}$($X_{k},~X_{k-1},~\ldots,~X_{k-w}$)\; } $k=k+1$ }%chktex 36
\end{algorithm2e}
The basic algorithm for the RRE-accelerated iterative solution of~\eqref{sylmateq} is summarized in Algorithm~\ref{alg:basic}. A slightly different formulation of the problem can be found in~\cite{JaMPR18} and, for the generalized Lyapunov case in, e.g.,~\cite{Dam08}.
It requires the solution of a standard Sylvester equation in every step. For small to moderately sized problems, this can be achieved by, e.g., the Bartels-Stewart algorithm~\cite{BarS72} which operates by first transforming $A$, $B$ to their (real) Schur forms
\[
  Q^{\TT}_A A Q_A = R_A,\quad Q_B^{\TT} B Q_B = R_B,
\]
where $Q_A\in\Rnn,~Q_B\in\Rmm$ are unitary
and $R_A\in\Rnn,~R_B\in\Rmm$ quasi upper triangular. Then the transformed Sylvester equations with the quasi-triangular matrices $R_A,~R_B$ are solved.

In the course of Algorithm~\ref{alg:basic}, the operator $\cL$ does not
change. Hence, for efficiency reasons, it is paramount to compute these Schur decompositions once and apply them to the whole multi-term equation before initiating the iteration:
\begin{align*}
  A&\leftarrow R_{A}=Q_{A}^{\TT}AQ_{A},\quad B\leftarrow R_{B}=Q^{\TT}_{B}BQ_{B},\quad Y\leftarrow Q_{A}^{\TT}YQ_{B},\quad X\leftarrow Q_{A}^{\TT}XQ_{B},\\
  N_{k}&\leftarrow Q_{A}^{\TT}N_{k}Q_{B},\quad H_{k}\leftarrow Q_{B}^{\TT}H_{k}Q_{B},\quad k=1,\ldots,\ell.
\end{align*}

The stationary iteration is then applied to the transformed equation with quasi-triangular $A,~B$, which drastically reduces the costs to solve the standard Sylvester equations in every step.
Afterwards, a back-transformation is carried out to recover the approximate solution as $X_k\leftarrow Q_{A}X_{k}Q_{B}^{\TT}$.

In case of generalized Sylvester operators $\cL(X)=AXC+MXB$ with
$M\in\Rnn,~C\in\Rmm$, QZ-decompositions have to be used: $AZ_{A}=R_{A}Q_{A}$,
$MZ_{A}=T_{A}Q_{A}$, $BZ_{B}=R_{B}Q_{B}$, $CZ_{B}=T_{B}Q_{B}$, see, e.g.~\cite[Chapter 6]{Koe21},~\cite{Koe-mepack-all-versions} for details.

%%%%%%%%%%%%%%%%%%%%%%%%%%%%%%%%%%%%%%%%%%%%%%%%%%%%%
\subsection{Incorporation of RRE}\label{ssec:RRE_smalldense}
The minimization problem for finding the updated \(X_{k}\) via RRE from the
sequence \(\{X_{i}\}_{i=1}^{w}\) generated by Algorithm~\ref{alg:basic} is to find the coefficient vector \(\widehat{\gamma}\) satisfying
\begin{align}\label{eq:rre-dense-mat}
  \widehat{\gamma} = \underset{\gamma\in\Gamma}{\argmin}
  \left\| \sum_{i=1}^{w} \gamma_i (X_{i+1}-X_i) \right\|,
\end{align}
recalling that  \(\Gamma:=\left\lbrace
  \gamma\in\R^{w}|	\sum_{i=1}^{w} \gamma_i = 1 \right\rbrace\).
It is straightforward to handle these problems by vectorization
\begin{align*}
  \widehat{\gamma} = \underset{\gamma\in\Gamma}{\argmin}
  \left\| \sum_{i=1}^{w} \gamma_i (\vecop{(X_{i+1})}-\vecop{(X_i)}) \right\|,
\end{align*}
which can be solved exactly as in~\eqref{eq:rre-minprob}.
This formulation leads to normal equations defined by a forward difference matrix $U\in\R^{nm\times w}$. Clearly, this strategy is only feasible for small to moderate dimensions $n,~m$ and for $w\ll \min{(n,m)}$.
%%%%%%%%%%%%%%%%%%%%%%%%%%%%%%%%%%%%%%%%%%
\subsection{Convergence considerations}\label{ssec:conv_small}
To analyze the convergence of the splitting-based fixed-point iteration~\eqref{eq:basicit} with RRE, we rewrite it in vector form.
Let $x=\vecop(X)$ and $y=\vecop(Y)$. Using standard properties of the Kronecker
product \(\otimes\), the application of the operator $\cL(X)=AX+XB$ to \(X\)
corresponds to the multiplication of \(x\) with the matrix $\cLvec=I_m\otimes A
+ B^{\TT}\otimes I_n$. Similarly, $\Pi(X)=\sum_{k}N_{k} X H_{k}$ corresponds to $\Pivec=\sum_{k} (H_{k}^{\TT}\otimes N_{k})$. Both matrices are of size $b\times b$, where $b:=mn$.
The fixed-point iteration is therefore equivalent to
\begin{align}\label{eq:basicit-vec}
  x_{k+1}= \mathcal{G} x_{k} + c,\quad \mathcal{G} := -\cLvec^{-1}\Pivec \in\R^{b\times b}, \quad c = -\cLvec^{-1}(y)\in\R^{b}.
\end{align}
Suppose the eigenvalues $\lambda_i$ of $\mathcal{G}$ are ordered by modulus, i.e.
\begin{equation*}
  |\lambda_1| \ge |\lambda_2| \ge \dots \ge |\lambda_{b}|.
\end{equation*}
Convergence of the unaccelerated iteration~\eqref{eq:basicit-vec} requires the spectral radius $\rho(\mathcal{G}) = |\lambda_1| <1$, in which case the standard convergence theory implies~\cite[sect.~4.2]{Saad03}
\begin{equation}\label{eq:basicit-ek}
  e_{k} = \mathcal{G}^{k} e_0, \quad e_k := x_k-x, \quad k\ge 0.
\end{equation}
Denote the $w$-term extrapolant (starting with $x_k$) and the associated
difference matrix as
\begin{equation*}
  \hx_{k,w} =   \sum_{i=1}^{w} \hgamma_i x_{k+i-1}, \quad
  U_{k,w} = [u_k, u_{k+1}, \ldots, u_{k+w-1}],
\end{equation*}
where the first-order difference vector $u_i$ is defined in~\eqref{eq:def-ui} and
\begin{equation*}
  \hgamma = \underset{\gamma\in\Gamma}{\argmin}
  \left\| U_{k,w}\cdot \gamma \right\|_2 =
  \underset{\gamma\in\Gamma}{\argmin} \left\|  \sum_{i=1}^{w} \gamma_i u_{k+i-1} \right\|_2.
\end{equation*}
Using the formulation of Sidi~\cite{sidi1986convergence}
the extrapolated iterate can be alternatively written as
\begin{equation}\label{eq:basicit-x-kw}
  \hx_{k,w} =   x_k + \sum_{i=1}^{w-1} \hq_i u_{k+i-1}= x_k + U_{k,w-1}\hq, \quad
  \hq = [\hq_1, \hq_2,\dots,\hq_{w-1}]^{\TT}\in\R^{w-1},
\end{equation}
where the coefficients $\hq_i$ satisfy
\begin{equation*}
  \hgamma_1 = 1-\hq_1, \quad \hgamma_{w}=\hq_{w-1}, \quad
  \hgamma_i = \hq_{i-1} - \hq_{i}, \quad i=2,\dots,w-1.
\end{equation*}
Let \(\Psi_{k,w} := [\psi_{k}, \psi_{k+1}, \ldots, \psi_{k+w-1}]\), where $\psi_i:= u_{i+1}-u_i$ are the vectors of second-order difference of the $x_i$.
Then, by construction $\hq$ solves the \emph{unconstrained} minimization problem
\begin{equation*}
  \hq = \underset{q\in\R^{w}}{\argmin} \| u_k + \Psi_{k,w-1}\cdot q  \|_2.
\end{equation*}

By using an asymptotic expansion of determinant expression for the extrapolated iterate and the corresponding error $\he_k := \hx_{k,w}-x$, Sidi~\cite[Thm.~3.1]{sidi1986convergence} shows, if $\mathcal{G}$ is diagonalizable, then, for some constant $C_\mathcal{G}>0$ independent of $k$,
\begin{equation}\label{eq:basicit-rre-rate}
  \normt{\hx_{k,w} - x} \le
  C_\mathcal{G}\left(
	\frac{\lambda_{w}}{\lambda_1}\right)^k \normt{x_{k+w} -x} \quad
  \text{as}\quad
  k\to \infty.
\end{equation}
This means the asymptotic reduction rate for~\eqref{eq:basicit-vec} with RRE becomes $|\lambda_{w}|$, giving significantly faster convergence than the non-accelerated rate $|\lambda_1|$ whenever $|\lambda_w|<|\lambda_1|<1$.
More importantly, even if the sequence $\{x_i\}$ is divergent, i.e., $|\lambda_{1}|\ge 1$, $\hx_{k,w} \to x$ as $k\to\infty$, provided that $|\lambda_{w}|<1$.
On the other hand, the bound~\eqref{eq:basicit-rre-rate} implies that no
significant convergence acceleration might take place when applying RRE if
$|\lambda_1| \approx |\lambda_2| \approx \dots \approx |\lambda_{w}|$. From a
practical point of view, finding a suitable or optimal value for $w$ is
difficult, since it requires knowledge of the full spectrum of $\cG$ which is not efficiently available.

From~\eqref{eq:basicit-ek}, we have
\begin{equation*}
  u_{k} = x_{k+1} - x_{k} = e_{k+1} - e_{k} = (\mathcal{G}-I) e_k = (\mathcal{G}-I)\mathcal{G}^k e_0.
\end{equation*}
Substituting this equation into~\eqref{eq:basicit-x-kw} and using~\eqref{eq:basicit-ek} gives, assuming that $\mathcal{G}$ is regular,
\begin{equation}\label{eq:basicit-error}
  \hx_{k,w} - x = e_k + \sum_{i=1}^{w-1} \hq_i \left(\mathcal{G}-I\right) \mathcal{G}^{i-1}e_k
  = \left(
	\mathcal{G}^{-w} + \sum_{i=1}^{w-1} \hq_i \left(\mathcal{G}^{i-w}-\mathcal{G}^{i-w-1}\right)
  \right)e_{k+w} =: h(\mathcal{G}) e_{k+w},
\end{equation}
where we have defined the rational function $h(\lambda) = \lambda^{-w} + \sum_{i=1}^{w-1} \hq_i (\lambda^{i-w}-\lambda^{i-w-1})$.
Comparing~\eqref{eq:basicit-error} with the asymptotic bound~\eqref{eq:basicit-rre-rate}, we know $h(\lambda)$ is almost annihilating at the leading eigenvalues $\lambda_{1},\dots,\lambda_{w-1}$ of $\mathcal{G}$ as $k$ tends to infinity~\cite{SmiFS87}, that is,
\begin{equation}\label{eq:basicit-h}
  \norm{h(\mathcal{G})} \in \cO\left(
	\Big(\frac{\lambda_{w}}{\lambda_1} \Big)^k
  \right) \quad \text{as} \quad k\to \infty.
\end{equation}

A strategy that can gain further acceleration via RRE is using the cycling mode, in which case the fixed-point iteration~\eqref{eq:basicit-vec} is initialized in each cycle using the extrapolant of the previous cycle.
It is beneficial for an efficient implementation, since only $w$ additional iterates need to be stored and the RRE coefficients are only computed after every $w$ steps of the stationary process.

The RRE extrapolated error in the cycling mode is discussed in~\cite{Sidi1992upper} for linear systems and in~\cite{sidi2020convergence} for nonlinear systems.
For the linear case, upper bounds in the cycling mode are given in terms of Jacobi polynomials on $\mathcal{G}$,
and the bound can be expressed using certain types of spectra of it~\cite[Thm.~7.1]{sidi2007pagerank},~\cite[sect. 7.3]{Sidi17}.
In the case where the cycling starts from $x_0$ and the window size $w_i$ at the
$i$th cycle is \emph{exactly} the degree of the minimal polynomial of the
Jacobian matrix at the fix point of the previous cycle (with respect to
$\hx^{(i-1)}_{0,w_{i-1}}$), then the quadratic convergence of the
sequence $\{ \hx^{(i)}_{0,w_{i}} \}$ follows heuristically from the
result for minimal polynomial extrapolation~\cite[Thm.~4.1]{stig1980periodic}.
Moreover, suppose $\alpha(I-\mathcal{G})$ has positive definite Hermitian part for some $\alpha\in\C$, $|\alpha|=1$. Then, for a fixed window size $w$, the sequence $\{ \hx^{(i)}_{0,w_{i}} \}$ is proven to converge at least geometrically~\cite[Thm.~7.2]{Sidi17}.

%%%%%%%%%%%%%%%%%%%%%%%%%%%%%%%%%%%%%%%%%%%%%%%%%%%%%%%%%%%%%%%%%%%%%%%%%%%%
\section{Analysis and algorithm for large-scale problems}
We now turn our attention to the case when the coefficients of~\eqref{sylmateq} are large matrices, in particular when $A,~B$ of the Sylvester operator~\eqref{sylmateq} are sparse matrices and the $N_k,~H_k$ in~\eqref{posoper} allow efficient matrix-vector products.
Often, the right-hand side is of low rank, $r:=\rank{Y}\ll\max(n,m)$, and is
given in (or can be brought in) factored form $Y=FTG^{\TT}$ with
$F\in\Rnr,~G\in\Rmr,~T\in\Rrr$, which we assume in the remainder. Then the
solution can also be expected to have a low numerical rank and, hence, can be well approximated by a low-rank factorization $X\approx Z_{L}DZ_{R}^{\TT}$ with $Z_{L}\in\R^{n\times z},~Z_{R}\in\R^{m\times z},~D\in\R^{z\times z}$, with $z\ll \min(n,m)$.
A large amount of numerical evidence backing up this expectation can be found in studies for solving general linear matrix equations~\eqref{genmateq}, see, e.g.,~\cite{Dam08,KreT11,BenB13,HaoS21,JaMPR18,ShaSS15,PalCS25,PalK2021}.
Theoretical investigations on this low numerical rank for certain special cases can be found, e.g, in~\cite{Gra04,KreT11,BenB13} and, in particular, in~\cite{JaMPR18} for the present situation~\eqref{sylmateq}.
It is straightforward to turn Algorithm~\ref{alg:basic} into a low-rank
algorithm that generates low-rank factors of $X$ by replacing the Sylvester
solver with a low-rank Sylvester solver and expressing the remaining operations in
terms of the low-rank factors, together with suitable rank truncation to keep
the memory demands low. For Lyapunov-plus-positive equations, this strategy, without extrapolation, has been investigated in, e.g.,~\cite{Dam08,ShaSS15}.

For large multi-term Sylvester equations, Algorithm~\ref{alg:nonstat} summarizes
this low-rank version of Algorithm~\ref{alg:basic}. In the following, we discuss
the main algorithmic ingredients of Algorithm~\ref{alg:nonstat}. The main focus
lies on the incorporation---and its consequences---of RRE\@.

\begin{algorithm2e}[tbp]
  \caption{%
    Non-stationary iteration with RRE for the generalized Sylvester equation~\eqref{sylmateq}.
  }%
  \label{alg:nonstat}
  \KwIn{%
    Matrices $A,~B,~Y=FTG^{\TT},~\{N_i\}_{i=1}^{\ell},~\{H_i\}_{i=1}^{\ell},~X_0$, tolerances $\tau_{\text{outer}}, \tau_{\text{inner}}, \tau_{\text{trunc}}$, rank compression routine $\cT(\cdot)$
  }
  \KwOut{%
    Approximate solution~$\widetilde X=Z_{L}DZ_{R}^{\TT}\approx X$.
  }

  $k=1$\;
  $F_0=F,~T_0=T,~G_0=G$\;
  \While{not converged}{
    Solve $AX_{k}+X_{k}B=-F_{k}T_{k}G_{k}^{\TT}$ approximately w.r.t. $\tau_{\text{inner}}$ for $Z_{L,k}D_{k}Z_{R,k}^{\TT}$\;
    % Y-\sum\limits_{i=1}^{\ell}N_i X_{k-1} H_i$\;
	$\tX_{k}=Z_{L,k}D_{k}Z^{\TT}_{R,k}\leftarrow\cT(Z_{L,k}D_{k}Z^{\TT}_{R,k},\tau_{\text{trunc}})$\label{alg:trunc_sol}\;
	\If{$k\geq w~\land~\mod(k,w)=0$}{
      Extrapolation step: $X_{k} \leftarrow \mathrm{RRE}_{w}(\tX_{k},~\tX_{k-1},~\ldots,~\tX_{k-w}$)\;
      Truncate again: $X_k=Z_{L,k}D_{k}Z^{\TT}_{R,k}\leftarrow\cT(X_k,\tau_{\text{trunc}})$\label{alg:trunc_rre}\;
	}
	Estimate residual norm\;
	\If{residual $< \tau_{\mathrm{outer}}$}{STOP}
	$F_{k}T_{k}G_{k}^{\TT}=\cT(FTG^{\TT}+\sum\limits_{i=1}^{\ell}N_i Z_{L,k}D_{k}Z_{R,k}^{\TT} H_{i},\tau_{\text{trunc}})$\label{alg:trunc_rhs}\;
	$k=k+1$
  }
\end{algorithm2e}

%%%%%%%%%%%%%%%%%%%%%%%%%%%%%%%%%%%%%%%%%%%%%%%%%%%%%%%
\subsection{Incorporation of RRE}\label{ssec:RRE_lowrank}
Clearly, simply using vectorization of the iterates, as in the small, dense setting, is not an efficient option for large-scale problems, especially those exhibiting low-rank structure.

For efficiently executing RRE for the low-rank matrix sequence generated by
Algorithm~\ref{alg:nonstat}, we use and adapt the machinery developed
in~\cite{DeKLM25}. For generalized Lyapunov equations, the iterates are
symmetric low-rank matrices, $X_{i}=Z_{i}D_{i}Z_{i}^{\TT}$, and~\cite[Alg.~4]{DeKLM25} can be used right away. Here, we generalize this to the
nonsymmetric situation for multi-term Sylvester equations. Let
$\tX_{i}=Z_{L,i}D_{i}Z_{R,i}^{\TT}$ be the truncated iterates of
Algorithm~\ref{alg:nonstat} with $\rank{\tX_{i}}=z_{i}\ll\min{(n,m)}$, i.e., the factors $Z_{L,i}$, $Z_{R,i}$ have $z_{i}$ columns.
Consider the sum of increments in the optimization~\eqref{eq:rre-dense-mat}:
\begin{align*}
  \Delta&=\sum_{i=1}^{w} \hgamma_{i} (\tX_{i+1}-\tX_{i})=\sum_{i=1}^{w} \hgamma_{i} (Z_{L,i+1}D_{i+1}Z_{R,i+1}^{\TT}-Z_{L,i}D_{i}Z_{R,i}^{\TT})\\
        &=\begin{bmatrix}
          Z_{L,1} & \cdots & Z_{L,w+1}
        \end{bmatrix}
          \begin{bmatrix}
            -\hgamma_{1} D_{1} \\
            & \cgamma_{1} D_{2} \\
            && \ddots \\
            &&& \cgamma_{w-1} D_{w} \\
            &&&& \hgamma_{w} D_{w+1}
          \end{bmatrix}
          \begin{bmatrix}
            Z_{R,1}^{\TT} \\
            \vdots \\
            Z_{R,w+1}^{\TT}
          \end{bmatrix},
\end{align*}
where \(\cgamma_{i}= \hgamma_{i}-\hgamma_{i+1}\), and thin QR-factorizations of the left and right factors
\begin{align}\label{eq:qr_both}
  \begin{split}
    \begin{bmatrix}
      Z_{L,1} & \cdots & Z_{L,w+1}
	\end{bmatrix}&=Q_{L}R_{L}=Q_{L}\begin{bmatrix}
      R_{L,1} & \cdots & R_{L,w+1}
	\end{bmatrix},\\
	Q_{L}&\in\R^{n\times \sum_{i}z_{i}},\quad Q^{\TT}_{L}Q_{L}=I,\quad R_{L,j}\in\R^{\sum_{i}z_{i}\times z_{j}},\\
    \begin{bmatrix}
      Z_{R,1} & \cdots & Z_{R,w+1}
	\end{bmatrix}&=Q_{R}R_{R}=Q_{R}\begin{bmatrix}
      R_{R,1} & \cdots & R_{R,w+1}
	\end{bmatrix},\\
	Q_{R}&\in\R^{m\times \sum_{i}z_{i}},\quad Q_{R}^{\TT}Q_{R}=I,\quad R_{R,j}\in\R^{\sum_{i}z_{i}\times z_{j}}.
  \end{split}
\end{align}
Consequently, in any unitarily invariant norm
$\|\Delta\|=\|Q_L^{\TT}\Delta Q_R\|$ (cf.~\cite[Lemma~3.1]{DeKLM25}), the minimization objective
in~\eqref{eq:rre-dense-mat} can be written as
\begin{align*}
  \|\Delta\|&=\left\|
              \begin{bmatrix}
                R_{L,1} & \cdots & R_{L,w+1}
              \end{bmatrix}
              \begin{bmatrix}
                -\hgamma_1 D_1 \\
                & \cgamma_1 D_2 \\
                && \ddots \\
                &&& \cgamma_{w-1} D_w \\
                &&&& \hgamma_w D_{w+1}
              \end{bmatrix}
              \begin{bmatrix}
                R_{R,1}^{\TT} \\
                \vdots \\
                R_{R,w+1}^{\TT}
              \end{bmatrix}\right\|\\
            &=\left\| \sum_{i=1}^{w} \hgamma_i (R_{L,i+1}D_{i+1}R_{R,i+1}^{\TT}-R_{L,i}D_{i}R_{R,i}^{\TT})\right\|.
\end{align*}
Since every summand is of size $z_i\times z_i$, and we expect $z_i\ll\min{(n,m)}$,
it now becomes viable to vectorize the terms $R_{L,i}D_{i}R_{R,i}^{\TT}$ and carry out RRE as in the dense situation:
\[
  \hgamma = \underset{\gamma\in\Gamma}{\argmin}
  \left\| \sum_{i=1}^{w+1} \gamma_i \left(\vecop{\left(R_{L,i+1}D_{i+1}R_{R,i+1}^{\TT}\right)}-\vecop{\left(R_{L,i}D_{i}R_{R,i}^{\TT}\right)}\right) \right\|.
\]
The matrix of forward differences is then $U_w\in\R^{(\sum_{i}z_{i})^2\times w}$.  The
resulting extrapolant
$\hX=\sum\limits_{i=1}^w\hgamma_{i}Z_{L,i}D_{i}Z_{R,i}^{\TT}$ can
be of rank up to $\sum\limits_{i=1}^{w}z_{i}$ which can be larger than the rank of the most recent iterate. It is beneficial to apply a rank truncation to $\hX$, which
can be done efficiently thanks to the already computed
QR-factorizations~\eqref{eq:qr_both} and an SVD of the inner factors
\[
  \hX=Q_{L}\left(R_{L}\diag{(\hgamma_{1} D_{1},\ldots,\hgamma_{w} D_{w})}R_{R}^{\TT}\right)Q_{R}^{\TT}=Q_{L}\left(U_{L}\Sigma U_{R}^{\TT}\right)Q_{R}^{\TT}\approx \underbrace{Q_{L}\tU_{L}}_{=\hZ_{L}}\underbrace{\tSigma}_{=\hD} \underbrace{\tU_{R}^{\TT}Q_{R}^{\TT}}_{=\hZ_{R}^{\TT}},
\]
where $\widetilde \Sigma=\diag{(\sigma_1,\ldots,\sigma_{\tilde z})}$, and $\tU_L$ and $\tU_R$ contain the  singular vectors associated with the dominant singular values $\sigma_i$ with respect to a prescribed relative truncation tolerance
\(0<\tau_{\text{trunc}}\ll 1\).

%%%%%%%%%%%%%%%%%%%%%%%%%%%%%%%%%%%%%%%%%%%%%%%%%%%%%%
\subsection{Iterative solution of the standard Sylvester equation}\label{ssec:inner_solve}
The arising standard Sylvester or Lyapunov equation~\eqref{eq:basicsylv} can be solved for a low-rank approximated solution
\[
  AX_{k}+X_{k}B=-F_{k}T_{k}G_{k}^{\TT},\quad X_{k}\approx Z_{L,k}D_{k}Z_{R,k}^{\TT}.
\]
Here, $F_{k}T_{k}G_{k}^{\TT}$ is a low-rank representation of $FTG^{\TT}+\Pi(X_{k-1})$, which is discussed in detail in~Section~\ref{ssec:compress} (cf.\ equation~\eqref{rhs_lowrank}).
There are different solvers available for this task, e.g., low-rank ADI methods~\cite{li2002low,benner2014self} or extended or rational Krylov subspace methods~\cite{CasR23,Sim07,DruS11}, see also the surveys~\cite{BenS13, Sim16a} w.r.t.\ the Lyapunov case.
We typically refer to the iterations for the inner matrix equation~\eqref{eq:basicsylv} as inner iteration, whereas Algorithm~\ref{alg:nonstat} is called the outer iteration.

Solving~\eqref{eq:basicsylv} for varying accuracies, e.g.,
\[
  \|S_k\|:=\|AX_{k}+X_{k}B+F_{k}T_{k}G_{k}^{\TT}\|\leq \tau_{k,\text{inner}}
\]
during the outer fixed point process, makes Algorithm~\ref{alg:nonstat} a non-stationary iteration. The results from~\cite{ShaSS15} for generalized Lyapunov equations can be easily adapted for generalized Sylvester equations (see below).
The main ingredient is to enforce the inner residual norm to be proportional to the outer residual norm:
\[
  \|S_k\|\leq\eta\|\cL(X_k)+\Pi(X_k)+FTG^{\TT}\|,\quad \eta>0.
\]
%%%%%%%%%%%%%%%%%%%%%%%%%%%%%%%%%%%%%%%%%%%%%%%%%%%%%%%%%
\subsubsection{Convergence considerations}
To analyze the convergence of the RRE-accelerated non-stationary fix-point iteration presented as Algorithm~\ref{alg:nonstat}, we again recast the iteration in vector form, as in Section~\ref{ssec:RRE_smalldense}.
The key difference here is that, to obtain the next iterate, each fixed-point iteration is not solved exactly, but to a prescribed accuracy.
Consider the non-stationary fixed-point iteration
\begin{equation}\label{eq:inexact-it-vec}
  x_{k}= \mathcal{G} x_{k-1} + c + \delta_k,\quad \delta_k := \cLvec^{-1}s_k \in\R^{b},
\end{equation}
where $s_k:= \cLvec x_k + \Pivec x_{k-1} + y$ denotes the inner residual vector, with $c$ and $\mathcal{G} = -\cLvec^{-1}\Pivec$ defined the same as in~\eqref{eq:basicit-vec}.
Using~\eqref{eq:inexact-it-vec}, the fact that \(x\) is a fixed-point
of~\eqref{eq:basicit-vec} and that \(x_{k-1} = e_{k-1}+x\), the error arising in this non-stationary iteration satisfies the recurrence
\begin{equation}\label{eq:inexact-it-ek}
  e_k := x_k - x = \mathcal{G} e_{k-1} + (\mathcal{G}x + c - x) + \delta_k
  = \mathcal{G} e_{k-1}  + \delta_k,
\end{equation}
from which it follows
\begin{equation}\label{eq:inexact-it-ekj}
  e_{k+\ell} = \mathcal{G}^\ell e_{k}  + \sum_{j=1}^{\ell} \mathcal{G}^{\ell-j}\delta_{k+j}.
\end{equation}
Define the vectorized residual of the equation~\eqref{sylmateq} as
\[
  \xi_k := (\cLvec + \Pivec)x_k + y = (\cLvec + \Pivec)e_k.
\]
Now, if we impose $\norm{s_k} \le \eta \norm{\xi_{k-1}}$ for some $\eta>0$, we get
\begin{equation*}
  \norm{\delta_k}  \le \norm{\cLvec^{-1}} \norm{s_k}
  \le  \eta \norm{\cLvec^{-1}} \norm{\cLvec + \Pivec} \norm{e_{k-1}},
\end{equation*}
and, from~\eqref{eq:inexact-it-ek}, the immediate norm bound
\begin{equation}\label{eq:inexact-it-norm-bound}
  \norm{e_k} \le  \norm{\mathcal{G}} \norm{e_{k-1}}  + \norm{\delta_k}
  \le (\norm{\mathcal{G}}  + \eta \norm{\cLvec^{-1}} \norm{\cLvec + \Pivec}) \norm{e_{k-1}}.
\end{equation}
The bound~\eqref{eq:inexact-it-norm-bound} associated with the generalized Lyapunov case is used to prove the result in~\cite[Thm.~1 (i)]{ShaSS15}, which says that the non-stationary iteration~\eqref{eq:inexact-it-vec} without RRE converges if
\[
  \norm{s_k} \le \eta \norm{\xi_{k-1}} \quad \text{and} \quad \|\mathcal{G}\|+\eta\|\cLvec^{-1}\|\|\cLvec+\Pivec\|<1.
\]
This is a rather restrictive condition that is stronger than requiring solely $\norm{\mathcal{G}}<1$ due to the inexactness of the solver in each fixed-point iteration.

We now consider the effect of RRE on the iteration~\eqref{eq:inexact-it-vec}.
With the first-order difference vectors defined as in Section~\ref{ssec:RRE_smalldense} by $u_{k} = x_{k+1} - x_{k}$, the $w$-term extrapolant $\widetilde x_{k,w}$ satisfies the same expression as~\eqref{eq:basicit-x-kw}.
We now have, from~\eqref{eq:inexact-it-ek},
\begin{equation*}
  u_k = e_{k+1} - e_k = (\mathcal{G}-I)e_k + \delta_{k+1},
\end{equation*}
and therefore, by~\eqref{eq:basicit-error} the $w$-term extrapolated error satisfies
\begin{equation*}
  \hx_{k,w} -x =  e_k + \sum_{i=1}^{w-1} q_i (\mathcal{G}-I)e_{k+i-1}
  + \sum_{i=1}^{w-1} q_i \delta_{k+i}.
\end{equation*}
From~\eqref{eq:inexact-it-ekj}, for \(1\le i\le w-1\), we have
\[
  e_k = \mathcal{G}^{-w} e_{k+w} -  \sum_{j=1}^{w} \mathcal{G}^{-j}\delta_{k+j}, \qquad
  e_{k+i-1} = \mathcal{G}^{-w+i-1} e_{k+w} -  \sum_{j=1}^{w-i+1} \mathcal{G}^{-j}\delta_{k+i-1+j}.
\]
By substituting these  identities and rearranging, we arrive at
\begin{align}\label{eq:inexact-it-x-kw}
  \hx_{k,w} -x =&
                  h(\mathcal{G}) e_{k+w} - \sum_{i=1}^{w}\mathcal{G}^{-i}\delta_{k+i} +
                  \sum_{i=1}^{w-1} q_i \Big(
                  \delta_{k+i}  - (\mathcal{G}-I)\sum_{j=1}^{w-i+1} \mathcal{G}^{-j}\delta_{k+i-1+j}
                  \Big) \nonumber \\
  =&
	 h(\mathcal{G}) e_{k+w} +  \sum_{i=1}^{w} \Big(
	 (\mathcal{G}^{-i} + q_i I)\delta_{k+i} -
	 q_i(\mathcal{G}-I) \sum_{j=1}^{w-i+1} \mathcal{G}^{-j}\delta_{k+i-1+j}
	 \Big) \nonumber \\
  =& h(\mathcal{G}) e_{k+w} + \sum_{i=1}^{w}\sum_{j=1}^{w-i+1} C_{ij}\delta_{k+i+j-1},
\end{align}
where the rational function $h$ is defined in~\eqref{eq:basicit-error}, and we have introduced $q_w=0$ (in addition to those $q_i$ in~\eqref{eq:basicit-x-kw}) and the coefficient matrices
\begin{equation*}
  C_{ij} \equiv C_{ij}( \mathcal{G},q_i) = \begin{cases}
    \mathcal{G}^{-i} + q_i  \mathcal{G}^{-1},& j = 1, \\
    -q_i( \mathcal{G}-I) \mathcal{G}^{-j},& j = 2,\dots,w-i+1.
  \end{cases}
\end{equation*}
Taking norms on both sides of~\eqref{eq:inexact-it-x-kw} and using~\eqref{eq:basicit-h} gives
\begin{align}\label{eq:inexact-it-x-kw-norm}
  \norm{\hx_{k,w} -x } &\le
                         \norm{h(\mathcal{G})} \norm{e_{k+w}} + \sum_{i=1}^{w}\sum_{j=1}^{w-i+1} \norm{C_{ij}}\norm{\delta_{k+i+j-1}} \nonumber \\
                       &= \xi \norm{e_{k+w}} + \sum_{i=1}^{w}\sum_{j=1}^{w-i+1} \norm{C_{ij}}\norm{\delta_{k+i+j-1}} \quad
                         \text{as} \quad k \to \infty,
\end{align}
where $\xi\in \cO\left(
  \Big(\frac{\lambda_{w}}{\lambda_1} \Big)^k\right)$.
This bound shows that the $w$-term extrapolated error of the non-stationary iteration~\eqref{eq:inexact-it-vec} annihilates the contributions from the first $w-1$ dominant eigenvalues of $\mathcal{G}$, as in the stationary case, but it also contains the propagated error from the inexact solves (in the summation term).
The bound~\eqref{eq:inexact-it-x-kw-norm} recovers Sidi's asymptotic bound~\cite[Thm.~3.1]{sidi1986convergence} when the inner solver used in each iteration is exact, that is, $\delta_k\equiv0$ for all $k$. In principle, it shows that RRE can accelerate the convergence of~\eqref{eq:inexact-it-vec} as in the stationary case if the inner solver is \emph{sufficiently} accurate.

The $\delta_k$ terms in~\eqref{eq:inexact-it-x-kw-norm} can be easily translated into the error vectors $e_{k-1}$ under assumptions that link the inner residual vector $s_k$ and the full residual of the equation $\xi_{k-1}$. For example, if we impose~\cite[Thm.~1 (ii)]{ShaSS15}
\begin{equation*}
  \norm{\mathcal{G}} <1 \quad \text{and} \quad
  \norm{\cLvec^{-1}s_k} \le \eta \norm{\cLvec^{-1} \xi_{k-1}},\; \eta>0,
\end{equation*}
then we have
\begin{equation*}
  \norm{\delta_k} \le \eta \norm{(I - \mathcal{G})e_{k-1}} \le
  \eta \norm{I - \mathcal{G}} \norm{e_{k-1}},
\end{equation*}
where $\eta \norm{I - \mathcal{G} }<1$ means the inexactness of the solver in the $k$th iteration is not larger than the error in the previous iterate.

Finally, note that the above analysis also accommodates the case where the right-hand side vector $c=-\cLvec^{-1}(\vecop{Y})$ is only accessible inexactly---for example, when the matrix $Y$ is obtained via a low-rank truncation~\cite{ShaSS15}.

%%%%%%%%%%%%%%%%%%%%%%%%%%%%%%%%%%%%%%%%%%%%%%%%%%%%%%%%%
\subsubsection{Further considerations w.r.t.\ the inner solves}\label{sssec:inner_solve}
In our experiments, as the inner low-rank solver for~\eqref{eq:basicsylv}, we either use low-rank ADI or extended Krylov subspace methods.
The complexity per iteration step of both methods depends on the rank $r_k$ of the right-hand side $F_{k}T_{k}G_{k}^{\TT}$:
in low-rank ADI methods, the main effort per iteration step lies in the solution of (shifted) linear systems defined by $A+\beta_{j}I_{n},~B+\alpha_{j}I_{m}$ and right-hand sides with $r_{k}$ columns each. The shift parameters $\alpha_{j},\beta_{j}$ required in ADI methods can be computed before or automatically during the course of the iteration.

Extended (block) Krylov subspace methods also require solving linear systems with $r_k$ columns in the right-hand sides, albeit with constant coefficient matrices $A,~B$. Hence, matrix factorizations of $A,~B$ only need to be computed once and can be used throughout the whole non-stationary iteration. However, an underlying orthonormal basis of the extended block Krylov subspaces needs to be constructed at a cost of roughly $\cO(n(\mathrm{it}_{EKSM}\cdot r_k)^2)$. After $\mathrm{it}_{EKSM}$ steps of EKSM, the underlying Galerkin projection approach requires solving a smaller, dense version of~\eqref{eq:basicsylv} of dimension $\mathrm{it}_{EKSM}\cdot r_k$, which comes at a numerical cost in $\cO((\mathrm{it}_{EKSM}\cdot r_k)^3)$.
%%%%%%%%%%%%%%%%%%%%%%%%%
\subsubsection{Right-hand side separation}
In~\cite[Section~4.2]{ShaSS15}, EKSM was used exclusively as inner solver and the authors proposed to
reduce the numerical cost by separating $F_{k}T_{k}G_{k}^{\TT}$ into a sum of $r_k$ rank-one terms and, due to linearity, one can then solve the $r_k$ Sylvester equations independently. The iterate $X_k$ is then the sum of all individual low-rank solutions.

Here, we use a slightly more general separation into $1\leq p\leq r_k$ parts with ranks $1\leq r_k^{(p)}\leq r_k$:
\begin{align}\label{rhs_sep}
  \begin{split}
    F_{k}&=[F_{k}^{(1)},\ldots,F_{k}^{(p)}],\quad T_{k}=\diag{(T_{k}^{(1)},\ldots,T_{k}^{(p)})},\quad G_{k}=[G_k^{(1)},\ldots,G_{k}^{(p)}],\\
         & F^{(j)}_{k}\in\R^{n\times r_{k}^{(j)}},~G^{(j)}_{k}\in\R^{m\times r_{k}^{(j)}},\quad T^{(j)}_{k}\in\R^{r_k^{(j)}\times r_{k}^{(j)}},\quad \sum\limits_{j=1}^{p}r_{k}^{(j)}=r_k,\\
    F_{k}T_{k}G_{k}^{\TT}&=\sum\limits_{j=1}^{p}F^{(j)}_{k}T^{(j)}_{k}(G^{(j)}_{k})^{\TT}.
  \end{split}
\end{align}
We then solve the $p$ inner matrix equations
\begin{align}\label{inner_sylv_sep}
  \cA(X_{k}^{j})=-F^{(j)}_{k}T^{(j)}_{k}(G^{(j)}_{k})^{\TT},\quad j=1,\ldots,p,
\end{align}
each defined by a right-hand side of lower rank, potentially decreasing the overall computational cost. Finally, the iterate is set to $X_{k}\approx \sum\limits_{j=1}^{p}X_{k}^{(j)}$.

In our tests, this separation was observed to be more beneficial when EKSM is employed as the inner solver, because the numerical effort of EKSM scales cubically with the rank of the right-hand side. Low-rank ADI methods scale roughly linearly with $r_k$ (it mainly influences the cost for the arising shifted linear systems), so this separation strategy does not always lead to a run time reduction.
%%%%%%%%%%%%%%%%%%%%%%%%%%%%%%%%%%%%%%%%%%%%%%%%%%%%%%%%%
\subsection{Rank truncation}\label{ssec:compress}
In Algorithm~\ref{alg:nonstat}, rank truncation is carried out to prevent excessive storage consumption for certain matrices. This is denoted by $\cT(X,~\tau_{\text{trunc}})$, where $\tau_{\text{trunc}}$ is a truncation tolerance so that $\|X-\cT(X,~\tau_{\text{trunc}})\|\leq \tau_{\text{trunc}}$. In particular, these truncations are often necessary for the approximate solution of~\eqref{eq:basicsylv} (line~\ref{alg:trunc_sol}), for the right-hand side (line~\ref{alg:trunc_rhs}) and after RRE has been executed (line~\ref{alg:trunc_rre}), where the last instance was already discussed in Section~\ref{ssec:RRE_lowrank}.

We outline the main strategies here for the right-hand side of~\eqref{eq:basicsylv}. Consider at step $k$ a previous iterate $X_{k-1}=Z_{L,k-1}D_{k-1}Z_{R,k-1}$ of $\rank{X_{k-1}}=z_{k-1}$. Since
\begin{align}\begin{split}\label{rhs_lowrank}
  -Y-\Pi(X_{k-1})&=-FTG^{\TT}-\sum\limits_{i=1}^{\ell}N_i Z_{L,k-1}D_{k-1}Z_{R,k-1}^{\TT} H_i=F_{k}T_{k}G^{\TT}_{k},\\
                 &F_k:=[F,N_1 Z_{L,k-1},\ldots,N_{\ell} Z_{L,k-1}]\in\R^{n\times r_k},~T_{k}=\diag{(T,I_{\ell}\otimes D_{k-1})}^{r_k\times r_k},\\
                 &G_k:=[G,H_1^{\TT} Z_{R,k-1},\ldots,H_{\ell}^{\TT} Z_{R,k-1}]^{m\times r_k},
\end{split}
\end{align}
where $r_k=r + \ell z_{k-1}$, the right-hand side of~\eqref{eq:basicsylv} will be of rank up to $r_k$.
However, the numerical rank is usually much smaller than $r_k$ so that compression techniques are practical.
A common strategy is to compute thin QR-factorizations of $F_{k}=Q_{F}R_{F},~G_{k}=Q_{G}R_{G}$, followed by an SVD $U\Sigma V^{\TT}=R_{F}T_{k}R_{G}^{\TT}$, which allows one to adjust the truncation error by keeping singular vectors associated with singular values with
\[
  \sigma_{i}\geq \sigma_{1}\tau_{\text{trunc}}.
\]
If only $r$ singular vectors and values are kept,
then compressed low-rank factors are given by $\tF_{k}=Q_{F}[u_{1},\ldots,u_{r}]$, $\tT_{k}=\diag{(\sigma_{1},\ldots,\sigma_{r})}$, $\tG_{k}=Q_{G}[v_{1},\ldots,v_{r}]$.

As discussed in~\cite{PalCS25}, explicitly setting up the factors $F_{k},~G_{k}$ before the QR-factorization can be very memory consuming. Hence, the truncation is carried out one block at a time by first compressing each term $N_{i} Z_{L,k-1}D_{k-1}Z_{R,k-1}^{\TT}H_{i}$ individually, and subsequently compressing their sum:
\begin{align*}
  F_{k,i}T_{k,i}G_{k,i}^{\TT}&=\cT(N_{i} Z_{L,k-1}D_{k-1}Z_{R,k-1}^{\TT}H_{i},~\tau_{\text{trunc}}),\quad i=1,\ldots, \ell, \\
  F_{k}T_{k}G_{k}^{\TT}&=\cT([F,F_{k,1},\ldots,F_{k,\ell}]\diag{(T,T_{k,1},\ldots,T_{k,\ell})}[G,G_{k,1},\ldots,G_{k,\ell}]^{\TT},~\tau_{\text{trunc}}).
\end{align*}
The QR factorizations for the individual truncations can be reused for the truncation of the final sum. Alternative approaches using randomized approximate QR factorizations~\cite{HalMT11} can be used here, as discussed in~\cite{PalCS25}.

The truncation tolerance $\tau_{\text{trunc}}$ can be dynamically adjusted in the course of the outer non-stationary iteration similarly as the residual tolerances $\tau_{\text{inner}}$, see~\cite[Cor.~1]{ShaSS15}. There exist further strategies for adjusting the tolerances $\tau_{\text{trunc}}$, $\tau_{\text{inner}}$, see, e.g.,~\cite{Jen23}.
It is also wise to define a maximal allowed rank or column dimension, $\max_{\textrm{col}}\ll \min{(n,m)}$, to limit the storage consumption.

\begin{remark}
  Note that if the separation strategy for the right-hand side from the previous
  section (cf.\ equations~\eqref{rhs_sep},~\eqref{inner_sylv_sep}) is used, one can also additionally truncate the approximate solutions of the individual inner Sylvester equations~\eqref{inner_sylv_sep}, $X_k^{j}=\cT(X_k^{j},~\tau_{\text{trunc}}\frac{r_k^{(j)}}{r_k})$, $j=1,\ldots,p$. While this leads to additional truncation costs at this stage, it can reduce the costs of the further steps, e.g., when truncating the combined $X_k$.
\end{remark}
%%%%%%%%%%%%%%%%%%%%%%%%%%%%%%%%%%%%%%%%%%%%%%%%%%%%%%%%%
\subsection{Estimation of the residual norm}\label{ssec:res_norm}
Computing $\|\cA(X_k)\|$ for large-scale matrix equations is a demanding task on its own. For the spectral norm, we compute an estimate of $\sigma_1(\cA(X_k))$ by an iterative method for large-scale SVDs. Such methods only require matrix-vector products with $\cA(X_k)$ and its transpose, which can be cheaply computed without explicitly forming the residual matrix:
\begin{align*}
  \cA(X_k)p%&=(AX_k+X_kB+\sum\limits_{k=1}^{\ell}N_kXH_k)p\\
  =&A(Z_{L,k}(D_{k}(Z_{R,k}^{\TT}p)))+Z_{L,k}(D_k(Z_{R,k}^{\TT}(Bp)))\\
            &+\sum\limits_{k=1}^{\ell}N_{k}(Z_{L,k}(D_{k}(Z_{R,k}^{\TT}(H_{k}p))))+F(T(G^{\TT}p)),\quad p\in\Rm,
\end{align*}
and analogous for $\cA(X_{k})^{\TT}q,~q\in\Rn$.

A few dominant singular values $\sigma_{j}$ of $\cA(X_{k})$ might be used to estimate the Frobenius norm via $\|\cA(X_{k})\|_{F}\approx \sum_{j}\sigma_{j}$, although it is in general not known in advance how many approximate singular values are needed for a sufficiently accurate norm estimate.

Alternatively, one can compute an SVD of $\cA(X_k)$ by exploiting the low-rank structure, see, e.g., the strategies in~\cite[Section~4.3]{ShaSS15}. Similar block-by-block or randomized approaches discussed in Section~\ref{ssec:compress} can be used here as well.

%%%%%%%%%%%%%%%%%%%%%%%%%%
\section{Numerical examples}\label{sct:numerical_examples}
\rowcolors{2}{black!10}{white}
\renewcommand{\arraystretch}{1.2}
\begin{table}[t]
  \centering
  \input{dense_ex1}
  \caption{Result of dense algorithm with and without cycling RRE for the example in~\eqref{ex:1} with different values for $\beta,~\ell,~w$.}%
  \label{tab:dense1}
\end{table}

\begin{figure}[ht]
  \centering
  \input{multiterm_1}~\input{multiterm_1a}~\input{multiterm_1b}

  \ref{legfourone}%chktex 2

  \input{multiterm_2}~\input{multiterm_2a}~\input{multiterm_3}

  \caption{Residual norm history of dense algorithm with and without cycling RRE (RRE$_w$) for different values of the weighting parameter $\beta$, the number \(\ell\) of terms in \(\Pi\), and RRE window size \(w\).}%
  \label{fig:dense1}
\end{figure}
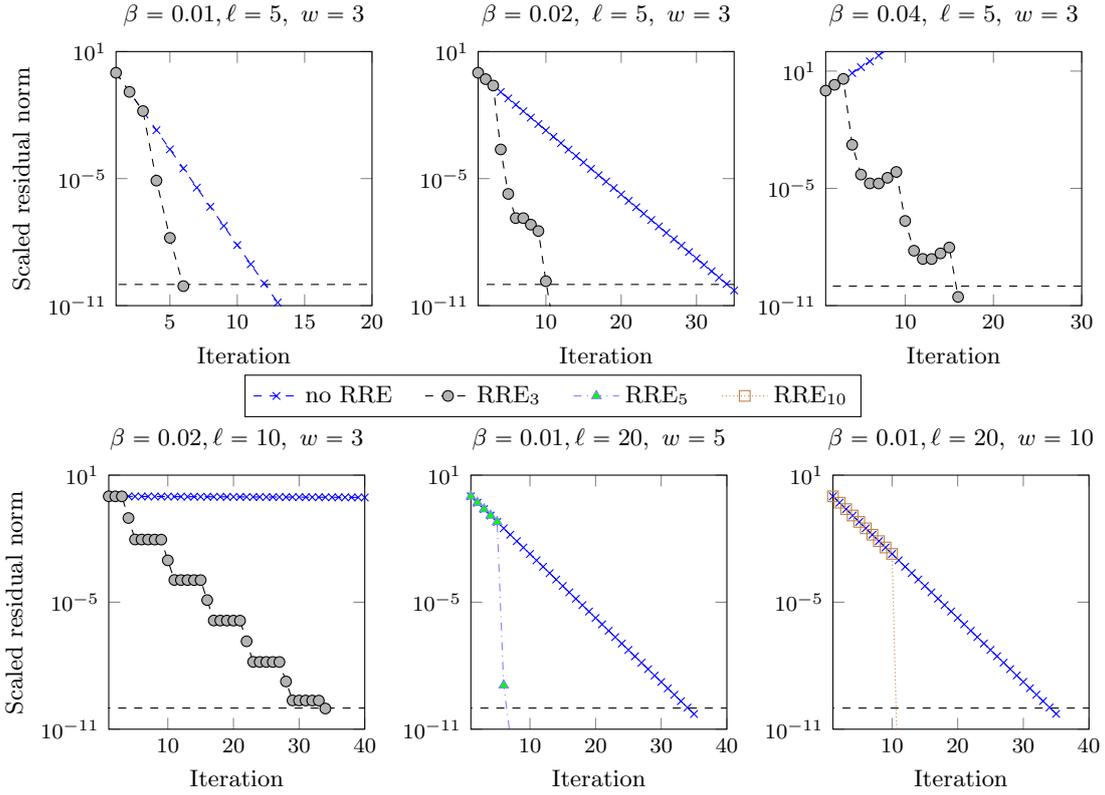

\begin{figure}[ht]
  \centering
  \input{gLyap_SSS}
  \small
  \begin{tabular}[t]{l|c|rrrr}
    \rowcolor{white}
    setting&inner&iter&$\rank{X}$&res&time [s]\\
    \hline
    no RRE~\ref{fig:advdiff:norre} &\cellcolor{white} &34 & 265 & \num{5.7917e-11} & \num{90.9393} \\
    RRE$_3$~\ref{fig:advdiff:rre3}& ADI &  25 & 265 & \num{3.6824e-11} & \num{72.5029} \\
    RRE$_5$~\ref{fig:advdiff:rre5}& \cellcolor{white} &\textbf{21} & 261 & \num{3.6783e-11} & \textbf{\num{56.1427}} \\
    \hline
    no RRE~\ref{fig:advdiff:norre}& &34 & 265 & \num{5.7294e-11} & \num{160.7389} \\
    RRE$_3$~\ref{fig:advdiff:rre3}&\cellcolor{white} EKSM & 25  &     267   &    \num{2.9629e-11}  &  \num{133.2300}\\
    RRE$_5$~\ref{fig:advdiff:rre5}& & \textbf{20} & \textbf{247} & \num{1.8299e-11} & \num{98.6421} \\
    \hline
  \end{tabular}
  \caption{Results for the advection diffusion example. Residual norm history of low-rank non-stationary iteration algorithm without and with cycling RRE using window size $w$ (RRE$_w$) and tabular summary.}%
  \label{fig:advdiff}
\end{figure}
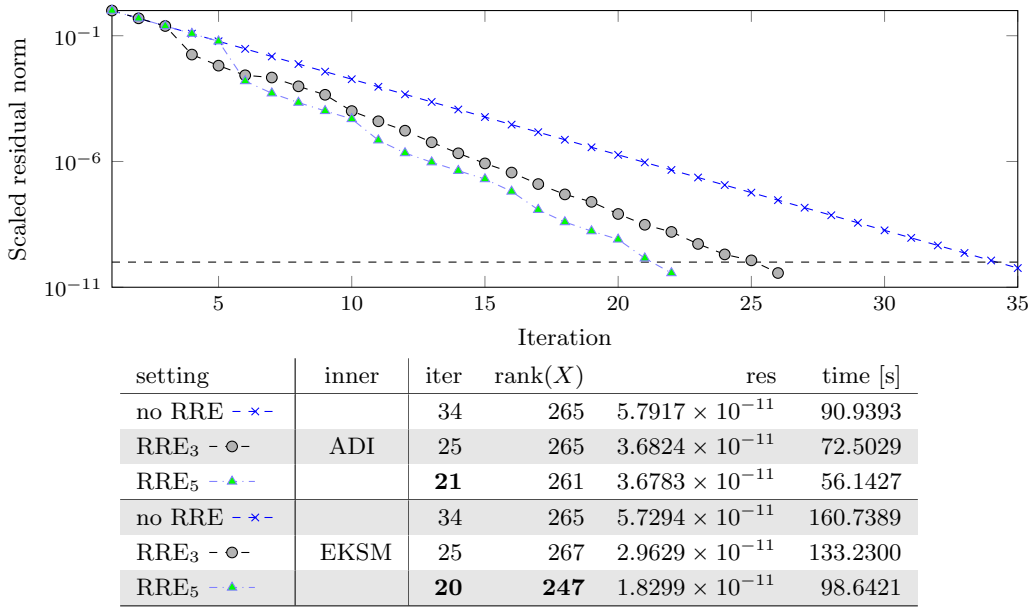

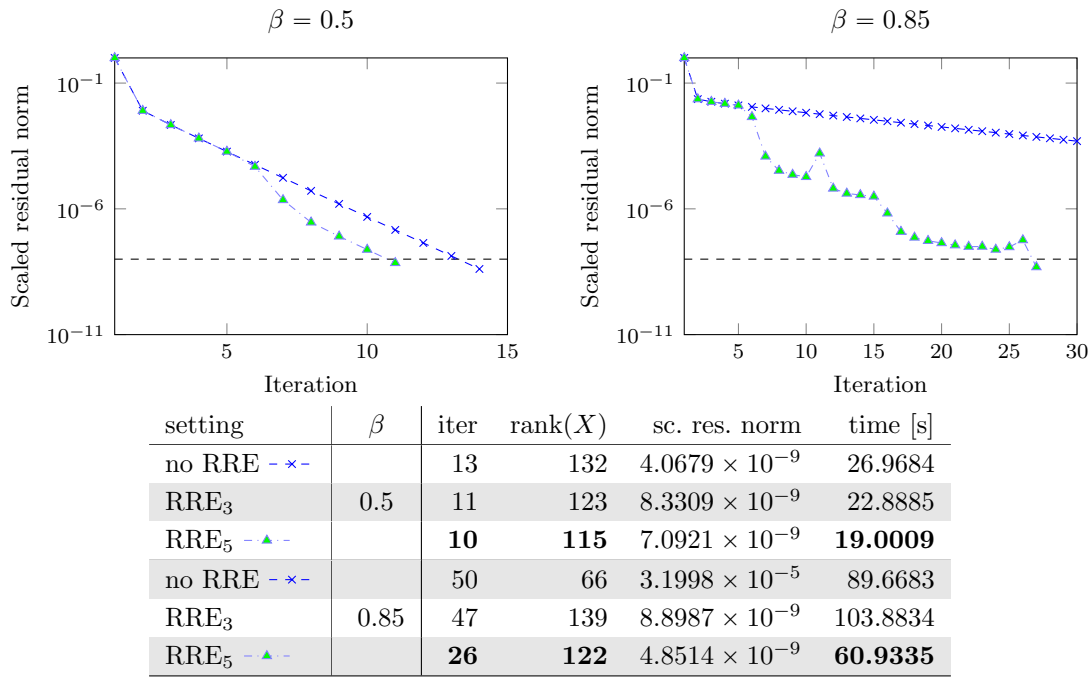
\begin{figure}[t]
  \centering
  \input{gLyap_RC}
  \begin{tabular}[t]{l|c|rrrr}
    \rowcolor{white}
    setting&$\beta$&iter&$\rank{X}$&sc.\ res.\ norm&time [s]\\
    \hline
    no RRE~\ref{fig:RCplot:norre}& \cellcolor{white}  &13 & 132 & \num{4.0679e-09} & 26.9684 \\
    RRE$_3$ & $0.5$  &  11 & 123 & \num{8.3309e-09} & 22.8885 \\
    RRE$_5$~\ref{fig:RCplot:rre5}& \cellcolor{white} &\textbf{10} & \textbf{115} & \num{7.0921e-09} & \textbf{19.0009} \\
    \hline
    no RRE~\ref{fig:RCplot:norre_2}& &50 & 66 & \num{3.1998e-05} & 89.6683 \\
    RRE$_3$ &\cellcolor{white} $0.85$ & 47 & 139 & \num{8.8987e-09} & 103.8834 \\
    RRE$_5$~\ref{fig:RCplot:rre5_2}&   & \textbf{26} & \textbf{122} & \num{4.8514e-09} & \textbf{60.9335} \\
    \hline
  \end{tabular}
  \caption{Results for the nonlinear circuit example for different $\beta$ values. Residual norm history of low-rank non-stationary iteration algorithm (for some selected cases) and tabular summary.}%
  \label{fig:circuit}
\end{figure}

\begin{figure}[t]
  \centering
  \input{gSylv_SSS}
  \begin{tabular}[t]{l|rrrr}
    \rowcolor{white}
    &iter&$\rank{X}$&res&time [s]\\
    \hline
    no RRE~\ref{fig:rreplot:norre}&32 & 230 & \num{5.2195e-11} & 96.0563 \\
    RRE$_3$~\ref{fig:rreplot:rre3}&  23  & 208  &   \num{5.8483e-11} &   72.1430 \\
    RRE$_5$~\ref{fig:rreplot:rre5}&   \textbf{20} & \textbf{208} & \num{3.0034e-11} & \textbf{62.8994} \\
    RRE$_8$~\ref{fig:rreplot:rre8}&    22 & 222 & \num{7.5414e-11} & 66.6835 \\
    \hline
  \end{tabular}
  \caption{Results for the multiterm Sylvester example for different $w$ values. Residual norm history of low-rank non-stationary iteration algorithm and tabular summary.}%
  \label{fig:multiterm_sylv}
\end{figure}

The experiments were run using the 64-bit Windows version of \matlab\ 2023a on a desktop computer equipped with an \amd~\textsf{Ryzen}~7~5800X\ 8-core processor at 3.80GHz and with 64GiB of RAM\@.
The goal in all experiments is to achieve $\|\cA(X_{i})\|_{2}\leq \epsilon\|FTG^{\TT}\|_{2}$, for some $\epsilon>0$, within $k_{\max}=50$ iteration steps of the (non-)stationary iteration. %chktex 36
In the dense version, Algorithm~\ref{alg:basic}, the preliminary Schur decompositions are computed using the \texttt{schur} command and the occurring standard Sylvester or Lyapunov equations are solved by the built-in routine~\texttt{lyap} of \matlab.
The low-rank iteration, Algorithm~\ref{alg:nonstat}, is for generalized Lyapunov equations based on extensions of the \mess{} package~\cite{SaaKB21-mmess-3.1}. The occurring standard Lyapunov equations are either handled by the extended Krylov subspace method or low-rank ADI iteration, both, provided by \mess. The ADI shift parameters are generated automatically during the first run of the LR-ADI method at step $k=1$ of Algorithm~\ref{alg:nonstat} and reused in all later runs.

The low-rank stationary iteration for multi-term Sylvester equations uses its own implementation, where the inner standard Sylvester equations are solved by the low-rank ADI iteration~\cite{BenK14,BenLT09} using the implementation available at~\cite{Kuersch_Sylvester_ADI}. The required ADI shifts are produced once a priori by the heuristic from~\cite{BenLT09} and used in all runs of LR-ADI\@.

If not stated otherwise, the low-rank iteration uses a dynamic inner solve tolerance with $\eta=10^{-3}$. Motivated by~\cite{ShaSS15} and the provided implementation, the truncation threshold is set in the same way.
The spectral norm of the residual matrix $\cA(X_k)$ is estimated using the \texttt{svds} command from \matlab.
\subsection{Code and data availability}
All codes and data required to reproduce and verify the numerical experiments reported in this section have been made available in~\cite{KueBS26zenodo}. The authors would like to thank Martin K{\"{o}}hler~\orcidlink{0000-0003-2338-9904} and Jonas Schulze~\orcidlink{0000-0002-2086-7686} for their careful reviewing of the codes and their verification of the experiments.%chktex 8
We also thank Tobias Breiten~\orcidlink{0000-0002-9815-4897} for providing codes to generate the matrices of the circuit example.%chktex 8

\subsection{An example with dense coefficients}
We generate random synthetic dense generalized Sylvester equations of dimension $n=500, m=300$, in the following way:
\begin{equation}\label{ex:1}
  \begin{aligned}
    A&=A_0-1.5\Re{\lambda_{\max}(A_0)}I_n,\qquad A_0=\mathrm{rand}(n),\\
    B&=B_0-1.5\Re{\lambda_{\max}(B_0)}I_m,\qquad B_0=\mathrm{rand}(m), \qquad Y=\mathrm{rand}{(n,m)},\\
    \Pi(X)&=\beta^2\sum\limits_{k=1}^{\ell}N_{k}XH_{k}, \qquad N_{k}=\mathrm{rand}(n),\qquad H_{k}=\mathrm{rand}(m).
  \end{aligned}
\end{equation}
Here, $\mathrm{rand}(\cdot)$ stands for a matrix with random entries drawn from a uniform distribution in $[0,1]$.
The shifting by $-1.5I$ is to generate stable  matrices $A,~B$. The parameter $\beta$ and the number $\ell$ of extra terms influence the contraction property of the splitting.

For several values of $\beta,~\ell$ and RRE window sizes $w$, the stationary iteration Algorithm~\ref{alg:basic} is executed until $\|\cA(X)\|_2\leq 10^{-10}\|Y\|_2$ or if at most $k_{\max}=50$ steps were carried out. We use the notation RRE$_w$ for RRE with windows size $w$. The results are summarized in Table~\ref{tab:dense1} and, for selected combinations of $\beta,~\ell,~w$, Figure~\ref{fig:dense1} shows the scaled residual norms $\|\cA(X)\|_2/\|Y\|_2$ against the iteration number.
We see that using RRE in all tested cases leads to a significant speed-up of the stationary iteration with respect to both the number of outer iterations and the run times. For some combinations of $\beta,~\ell$, the plain stationary iteration diverges because $\cL+\Pi$ is not a convergent splitting of $\cA$, while the accelerated iteration is still able to converge to the desired tolerance because the chosen window size $w$ is sufficiently large (cf.\ equation~\eqref{eq:basicit-rre-rate}). In some of the latter cases, the residual norm exhibits an oscillatory, yet decreasing, pattern; see Figure~\ref{fig:dense1} for $\beta=0.04,~\ell=5,~w=3$. Note that increasing $\beta$ even further will at some point lead to situations where RRE is not able to cure the divergence, even with very large window sizes $w$.

%%%%%%%%%%%%%%%%%%%%%%%%%%%%%%%%%%%%%%%%%%
\subsection{Large-scale examples}
Our examples for generalized Lyapunov equations ($\cL(X)=AX+XA^{\TT}$) are motivated by the aim to compute Gramians of bilinear control systems
\begin{align*}
  \dot{x}(t)&=A x(t)+\beta\sum\limits_{k=1}^\ell N_{k}u_{k}(t)x(t)+Bu(t).
\end{align*}
The parameter $\beta$, as well as $\ell$, influence the dominance of the Lyapunov operator $\cL$. The middle matrix on the right-hand side factorization is always set to $T=I_{r}$.

\subsubsection{Advection-diffusion model}\label{sssec:advdiff}
The model is a bilinear control system from a centered finite-difference discretization of $-\triangle \bx+\bx_y$ on $(0,1)^2$ with Robin-type boundary conditions $\bn^{\TT}\nabla \bx=\beta u(\bx-1)$ on the left- and right boundary but homogeneous Dirichlet boundary conditions otherwise. This leads to $\ell=2$ with varying sizes depending on the number of grid points. This setting is essentially the example \textsc{ADVDIFF} in~\cite{ShaSS15}, but we also use a different scaling parameter $\beta=0.8$, here, to make the equation more difficult to solve.

We consider a problem dimension $n=\num{22500}$ and test Algorithm~\ref{alg:nonstat} without RRE and with cycling RRE using different windows sizes $w$, the latter being denoted by RRE$_w$.
The goal is to achieve $\|\cA(X)\|_{2}\leq 10^{-10}\|Y\|_{2}$.
We also investigate the effect of the inner solver, namely EKSM and LR-ADI, in our experiments.
The right-hand sides for EKSM were separated into parts of at most rank $r_{k}^{(i)}=30$ (cf.\ Subsection~\ref{sssec:inner_solve}).
A separation into parts of smaller ranks led to higher computation times. When LR-ADI is the inner solver, no right-hand side separation was used.

The results are summarized in Figure~\ref{fig:advdiff}. The plot shows the history of the scaled residual norm $\|\cA(X)\|_2/\|Y\|_2$ generated by the non-stationary iteration without and with RRE$_3$, RRE$_5$ and LR-ADI as inner solver. We found that using EKSM as the inner solver almost does not change the plot.
It is obvious that RRE significantly speeds-up the iteration and reduces the  number of outer iteration steps.
The table in Figure~\ref{fig:advdiff} summarizes this trend and lists the rank, the scaled residual norm $\|\cA(X)\|_2/\|Y\|_2$ of the final iterate, and the computation times. Again, we see that RRE improves the iteration leading to fewer outer steps, smaller final ranks and reduced computation times. The best results are here obtained with RRE using $w=5$ and LR-ADI as inner solver. Using EKSM as inner solver does not affect the speed of the outer iteration, but generally led to higher overall computation times.

For the setting with RRE$_5$, we investigated where the major parts of the computation time is spent. If LR-ADI is used as inner solver, approximately \qty{60}{\percent} of the total time is spent for solving the shifted linear systems and roughly \qty{15}{\percent} are consumed by the rank truncations. The remaining percentages were attributable to other components, such as shift generation ($\approx \qty{8}{\percent}$) and norm computations (\qty{5}{\percent}).
With EKSM as inner solver, only roughly 8\% of the time is spent for solving linear systems (due to the reuse of the factorization of $A$) and  approximately \qty{12}{\percent} for rank truncations. The generation of the orthonormal basis including solving the projected Lyapunov equation consumed \qty{68}{\percent}. In both cases, the execution of RRE consumed less than \qty{1}{\percent} of the total run time.
%%%%%%%%%%%%%%%%%%%%
\subsubsection{Nonlinear circuit model}
This example stems from a nonlinear RC circuit model from~\cite{BenB13,morBaiS06} with $\ell=1$ extra terms and $r=1$.
We use the model with size $n=\num{22650}$ and different values $\beta=0.5$ or $\beta=0.85$. The goal is to achieve $\|\cA(X)\|_2\leq 10^{-8}\|Y\|_2$. Here, only LR-ADI is used as inner solver.
The results are summarized in Figure~\ref{fig:circuit}. In both cases, RRE improves the convergence speed of the non-stationary iteration, leading to reduced iteration numbers, ranks, and computation times, where RRE with $w=5$ gives the best result. Here, with $\beta=0.5$ the generalized Lyapunov equation is easier to solve by the non-stationary iteration, so that it converges sufficiently fast even without RRE and RRE allowing a modest speed-up.
The problem becomes more difficult with $\beta=0.85$, where the iteration without extrapolation converges very slow and fails to hit the target within $k_{\max}=50$ steps. Using RRE substantially improves the behavior and the goal is met in fewer than $k_{\max}$ steps. Again, a window size $w=5$ yields the best results.

%%%%%%%%%%%%%%%%%%%%
\subsubsection{A multi-term Sylvester equation}
We set up a multi-term Sylvester equation by using two different versions of the \textsc{ADVDIFF} example from Section~\ref{sssec:advdiff}.
The matrices $A,N_k,F$  and $B,H_k,G$ have dimension $n=\num{22500}$ and $m=\num{8100}$, respectively.
Here, $\ell=2$, $r=2$, and $\beta=0.8$. We execute the Sylvester version of Algorithm~\ref{alg:nonstat} without and with RRE using different window sizes $w$ and target accuracy $\|\cA(X)\|_2\leq 10^{-10}\|Y\|_2$. The results are summarized in Figure~\ref{fig:multiterm_sylv}. Once again, the RRE-accelerated non-stationary iteration is faster than without RRE\@. Increasing the window size $w$ improves the performance up to a certain point after which it slightly deteriorates (but it remains faster than without RRE). An optimal value for $w$ is difficult to predict since one usually has no knowledge about the dominant eigenvalues of the iteration map. Also note that a larger window size increases the storage needed for the required previous iterates.

%%%%%%%%%%%%%%%%%%%%%%%%%%%%%%%%%%%%%%%%%%%%%%%%%%%%%%%%%%%%%%%%%%%%%%%
\section{Conclusion and outlook}\label{sct:conclusion}
In this work we have discussed the RRE acceleration of iterative, fixed-point type, methods for multi-term Sylvester equations.
Both our theoretical analysis and the experimental validation showed that RRE can improve the performance of these methods, saving both memory and computing time.
For large-scale matrix equations, we extended the RRE framework from~\cite{DeKLM25} to sequences of nonsymmetric low-rank matrices. This was then incorporated into a non-stationary low-rank iteration for the generalized Sylvester equation. The machinery developed for the non-accelerated iteration for the symmetric case (generalized Lyapunov equations) in~\cite{ShaSS15}, such as inexact solutions of the occurring standard Sylvester equations and rank truncation with dynamically adjusted tolerances still work for the RRE accelerated iteration.

RRE uses $w$ previous iterates for executing the extrapolation, finding the ideal window size $w$ is a difficult task in practice and would be a potential further research topic. The theory predicts that RRE can in some cases generate convergent sequences even when the convergent splitting property is mildly violated~\cite{sidi1986convergence}. This was also observed in our experiments. Further research could investigate to what extent this behavior can be guaranteed.

In this work, we exclusively considered the situation where the operator $\cL(X)$ is a standard Sylvester or Lyapunov operator. However, the presented algorithms as well as the acceleration by RRE can be extended to other cases, provided efficient algorithms are available for solving matrix equations $\cL(X)=Y$. For example, if $\cL$ is the sum of a Sylvester operator and extra linear terms defined by low-rank coefficients, one can employ Sherman-Morrison-Woodbury approaches~\cite{Dam08,HaoS21} for solving with $\cL$.

%% file: dense_ex1.tex
\begin{tabular}{rrr | rrr | rrr}
  \rowcolor{white}
  &&&&no RRE&&&RRE(\(w\))&\\%chktex 36
  $\beta$ & $\ell$ & $w$ & iter & res & time & iter & res & time \\
  \hline
  0.01 & 5  & 3  & 12 & \num{1.4008e-11} & 1.0680 & 5 & \num{8.3324e-11} & 0.4703 \\
  0.02 & 5  & 3  & 34 & \num{5.2109e-11} & 2.4233 & 10 & \num{1.1564e-12} & 0.7954 \\
  0.04 & 5  & 3  & \multicolumn{3}{c|}{divergent} & 15 & \num{2.7845e-11} & 1.1773 \\
  0.02 & 10 & 3  & \multicolumn{3}{c|}{stagnation} & 33 & \num{9.4329e-11} & 2.9299 \\
  0.02 & 15 & 3  & \multicolumn{3}{c|}{divergent} & 15 & \num{8.6090e-11} & 1.7988 \\
  0.02 & 20 & 3  & \multicolumn{3}{c|}{divergent} & 16 & \num{6.3351e-12} & 2.3582 \\
  0.01 & 20 & 3  & 34 & \num{5.1206e-11} & 4.5641 & 9 & \num{3.3369e-11} & 1.4790 \\
  0.01 & 20 & 5  & 34 & \num{5.2335e-11} & 4.5019 & 6 & \num{4.5979e-12} & 1.0422 \\
  0.01 & 20 & 10 & 34 & \num{5.2925e-11} & 4.7518 & 10 & \num{1.8409e-14} & 1.6987 \\
  \hline
\end{tabular}

%% file: multiterm_1.tex
\begin{tikzpicture}
  \pgfplotstableread{multiterm_5_0.01_3.csv}\rredata
  \begin{semilogyaxis}[
    width  = 0.33\linewidth,
    height = 0.33\linewidth,
    ylabel = {Scaled residual norm},
    xlabel = Iteration,
    xmin   = 1,
    xmax   = 20,
    ymin   = 1e-11,
    ymax   = 10,
    % only marks,
    title  = {$\beta=0.01, \ell=5,~w=3$},
    title style = {font=\small},
    cycle list name = marks,
    ]
    \addplot+[] table[x expr = {\coordindex + 1}, y index = 0] {\rredata};
    \label{fig:mtone:norre}
    \addplot+[] table[x expr = {\coordindex + 1}, y index= 1] {\rredata};
    \label{fig:mtone:rre}
    \addplot[domain=0:50, dashed, thin, sharp plot, color=black] {10^ (-10)};
  \end{semilogyaxis}
\end{tikzpicture}%

%% file: multiterm_1a.tex
\begin{tikzpicture}
  \pgfplotstableread{multiterm_5_0.02_3.csv}\rredata
  \begin{semilogyaxis}[
    width  = 0.33\linewidth,
    height = 0.33\linewidth,
    xlabel = Iteration,
    xmin   = 1,
    ymin   = 1e-11,
    ymax   = 10,
    xmax   = 35,
    % only marks,
    title  = {$\beta=0.02,~\ell=5,~w=3$},
    title style={font=\small},
    cycle list name = marks,
    ]
    \addplot+[] table[x expr = {\coordindex + 1}, y index = 0] {\rredata};
    \label{fig:mtonea:norre}
    \addplot+[] table[x expr = {\coordindex + 1}, y index= 1] {\rredata};
    \label{fig:mtonea:rre}
    \addplot[domain=0:35, dashed, thin, sharp plot, color=black] {10^ (-10)};
  \end{semilogyaxis}
\end{tikzpicture}%

%% file: multiterm_1b.tex
\begin{tikzpicture}
  \pgfplotstableread{multiterm_5_0.04_3.csv}\rredata
  \begin{semilogyaxis}[
    width  = 0.33\linewidth,
    height = 0.33\linewidth,
    xlabel = Iteration,
    xmin   = 1,
    xmax   = 30,
    ymin   = 1e-11,
    ymax   = 100,
    % only marks,
    title  = {$\beta=0.04,~\ell=5,~w=3$},
    title style={font=\small},
    cycle list name = marks,
    % #2,
    %legend entries={{no RRE}, {RRE(\(3\))}, {RRE(\(5\))}, {RRE(\(10\))}},
    legend style={
      legend cell align=left,
      align = left,
      draw = black,
      fill = white,
      font = \small,
      legend columns=-1,
     legend to name=legfourone,
      /tikz/every even column/.append style={column sep=1em}
    }
    ]

    \addplot+[] table[x expr = {\coordindex + 1}, y index = 0] {\rredata};
    \addlegendentry{no RRE}
    \addplot+[] table[x expr = {\coordindex + 1}, y index= 1] {\rredata};
    \addlegendentry{RRE$_3$}
    \addplot+[mark phase=2] coordinates {(20,1e-10)}; % fake plot for extra legend entry
    \addlegendentry{RRE$_5$}
    \pgfplotsset{cycle list shift = 1}
    \addplot+[mark phase=2] coordinates {(20,1e-10)}; % fake plot for extra legend entry
    \addlegendentry{RRE$_{10}$}
    %\addlegendimage{\ref{fig:advdiff:rre5}}
    \addplot[domain=0:30, dashed, thin, sharp plot, color=black] {10^ (-10)};

  \end{semilogyaxis}
\end{tikzpicture}%

%% file: multiterm_2.tex
\begin{tikzpicture}
  \pgfplotstableread{multiterm_10_0.02_3.csv}\rredata
  % \createcolinterpolated\rredata
  % \pgfplotstableread{benner2020_q1_nonlinear_cycling.csv}\rredatacycling
  % \createcolinterpolated\rredatacycling
  % number of iterations for the cycling
  % \pgfplotstablegetrowsof{\rredatacycling}
  % \pgfmathsetmacro{\N}{\pgfplotsretval-1}
  \begin{semilogyaxis}[
    width  = 0.33\linewidth,
    height = 0.33\linewidth,
    ylabel = {Scaled residual norm},
    xlabel = Iteration,
    xmin   = 1,
    xmax   = 40,
    ymin   = 1e-11,
    ymax   = 10,
    % only marks,
    title  ={$\beta=0.02, \ell=10,~w=3$},
    title style={font=\small},
    cycle list name = marks,
    ]

    \addplot+[] table[x expr = {\coordindex + 1}, y index = 0] {\rredata};
    % \label{fig:rreplot:radi}
    \addplot+[] table[x expr = {\coordindex + 1}, y index= 1] {\rredata};
    % \label{fig:rreplot:noncycling}
    \addplot[domain=0:60, dashed, thin, sharp plot, color=black] {10^ (-10)};
  \end{semilogyaxis}
\end{tikzpicture}%

%% file: multiterm_2a.tex
\begin{tikzpicture}
  \pgfplotstableread{multiterm_20_0.01_5.csv}\rredata
   \begin{semilogyaxis}[
    width   = 0.33\linewidth,
    height  = 0.33\linewidth,
    xlabel  = Iteration,
    xmin    = 1,
    xmax    = 40,
    ymin    = 1e-11,
    ymax    = 10,
    title   = {$\beta=0.01, \ell=20,~w=5$},
    title style = {font=\small},
    cycle list name = marks,
    ]
    \addplot+[] table[x expr = {\coordindex + 1}, y index = 0] {\rredata};
    \pgfplotsset{cycle list shift = 1}
    \addplot+[] table[x expr = {\coordindex + 1}, y index = 1] {\rredata};
    \addplot[domain=0:60, dashed, thin, sharp plot, color=black] {10^ (-10)};
  \end{semilogyaxis}
\end{tikzpicture}%

%% file: multiterm_3.tex
\begin{tikzpicture}
  \pgfplotstableread{multiterm_20_0.01_10.csv}\rredata
  % \createcolinterpolated\rredata
  % \pgfplotstableread{benner2020_q1_nonlinear_cycling.csv}\rredatacycling
  % \createcolinterpolated\rredatacycling
  % number of iterations for the cycling
  % \pgfplotstablegetrowsof{\rredatacycling}
  % \pgfmathsetmacro{\N}{\pgfplotsretval-1}
  \begin{semilogyaxis}[
    width  = 0.33\linewidth,
    height = 0.33\linewidth,
    xlabel = Iteration,
    xmin   = 1,
    xmax   = 40,
    ymin   = 1e-11,
    ymax   = 10,
    % only marks,
    title  = {$\beta=0.01, \ell=20,~w=10$},
    title style={font=\small},
    cycle list name = marks,
    ]

    \addplot+[] table[x expr = {\coordindex + 1}, y index = 0] {\rredata};
    % \label{fig:rreplot:radi}
    \pgfplotsset{cycle list shift = 3}
    \addplot+[] table[x expr = {\coordindex + 1}, y index= 1] {\rredata};
    % \label{fig:rreplot:noncycling}
    \addplot[domain=0:60, dashed, thin, sharp plot, color=black] {10^ (-10)};
  \end{semilogyaxis}
\end{tikzpicture}%

%% file: gLyap_SSS.tex
\begin{tikzpicture}
  \pgfplotstableread{gLyap_SSS_22500_0.80_adi.csv}\rredata
  \begin{semilogyaxis}[
    width  = 0.9\linewidth,
    height = 0.35\linewidth,
    ylabel = {Scaled residual norm},
    xlabel = Iteration,
    xmin   = 1,
    xmax   = 35,
    ymin   = 1e-11,
    ymax   = 1,
    cycle list name = marks,
    ]
    \addplot+[] table[x expr = {\coordindex + 1}, y index = 0] {\rredata};
    \label{fig:advdiff:norre}
    \addplot+[] table[x expr = {\coordindex + 1}, y index= 1] {\rredata};
    \label{fig:advdiff:rre3}
    \addplot+[] table[x expr = {\coordindex + 1}, y index= 2] {\rredata};
    \label{fig:advdiff:rre5}
    \addplot[domain=1:50, dashed, thin, sharp plot, color=black] {10^ (-10)};
  \end{semilogyaxis}
\end{tikzpicture}%

%% file: gLyap_RC.tex
{
  \pgfplotsset{
    every axis/.style={
      width  = 0.45\linewidth,
      height = 0.35\linewidth,
      ylabel = {Scaled residual norm},
      xlabel = Iteration,
      xmin   = 1,
      xmax   = 30,
      ymin   = 1e-11,
      ymax   = 1,
      cycle list name = marks,
    },
  }
  \pgfplotstableread{gLyap_RC_22650_0.50_adi.csv}\rredata%
  \pgfplotstableread{gLyap_RC_22650_0.85_adi.csv}\rredatab%

  \begin{tikzpicture}
    \begin{semilogyaxis}[
      title = {\(\beta = 0.5\)},
      xmax=15
      ]
      \addplot+[] table[x expr = {\coordindex + 1}, y index = 0] {\rredata};
      \label{fig:RCplot:norre}
      \pgfplotsset{cycle list shift = 1}
      \addplot+[] table[x expr = {\coordindex + 1}, y index= 2] {\rredata};
      \label{fig:RCplot:rre5}
      \addplot[domain=1:50, dashed, thin, sharp plot, color=black] {10^ (-8)};
    \end{semilogyaxis}
  \end{tikzpicture}
  \hspace{1em}
  \begin{tikzpicture}
    \begin{semilogyaxis}[title = {\(\beta = 0.85\)}]
      \addplot+[] table[x expr = {\coordindex + 1}, y index = 0] {\rredatab};
      \label{fig:RCplot:norre_2}
      \pgfplotsset{cycle list shift = 1}
      \addplot+[] table[x expr = {\coordindex + 1}, y index= 2] {\rredatab};
      \label{fig:RCplot:rre5_2}
      \addplot[domain=1:50, dashed, thin, sharp plot, color=black] {10^ (-8)};
    \end{semilogyaxis}
  \end{tikzpicture}%
}

%% file: gSylv_SSS.tex
\begin{tikzpicture}
  \pgfplotstableread{gSylv_SSS_22500_8100_0.80.csv}\rredata
  \begin{semilogyaxis}[
    width  = 0.9\linewidth,
    height = 0.35\linewidth,
    ylabel = {Scaled residual norm},
    xlabel = Iteration,
    xmin   = 1,
    xmax   = 35,
    ymin   = 1e-11,
    ymax   = 1,
    cycle list name = marks,
    ]
    \addplot+[] table[x expr = {\coordindex + 1}, y index = 0] {\rredata};
    \label{fig:rreplot:norre}
    \addplot+[] table[x expr = {\coordindex + 1}, y index= 1] {\rredata};
    \label{fig:rreplot:rre3}
    \addplot+[] table[x expr = {\coordindex + 1}, y index= 2] {\rredata};
    \label{fig:rreplot:rre5}
    \addplot+[] table[x expr = {\coordindex + 1}, y index= 3] {\rredata};
    \label{fig:rreplot:rre8}
    \addplot[domain=1:50, dashed, thin, sharp plot, color=black] {10^ (-10)};
  \end{semilogyaxis}
\end{tikzpicture}%